\documentclass{amsart}

\usepackage[latin1]{inputenc}
\usepackage{amsmath}
\usepackage{amsthm}
\usepackage{amsfonts}
\usepackage{amssymb}
\usepackage[T1]{fontenc}
\usepackage[all]{xy}
\usepackage{graphicx}

\usepackage{calrsfs}
\usepackage{enumerate}
\newtheorem{p1}{P.}
\newtheorem{theo}{Th\'eor\`eme}

\newtheorem{coro}{Corollaire}
\newtheorem{exemple}{Exemple}

\begin{document}

\title[Pavages additifs]{Pavages additifs}
\author[J.-F. Marceau]{Jean-Fran\c cois Marceau}
\frenchspacing

\maketitle

Le présent article en est un d'introduction aux pavages additifs, faisant également des liens avec les pavages multiplicatifs. Nous ajoutons également deux nouveaux théorèmes pour les pavages additifs (voir théorème \ref{mien} et théorème \ref{mien 2}) et une conjecture sur les pavages multiplicatifs infinis.

\section{Introduction}
Au cours des dernières années, les avancées en matière d'algèbres amassées ont permis de révéler plusieurs liens avec divers domaines des mathématiques. Par exemple, un lien évident entre les frises de nombres générées à partir d'un carquois de type $A_n$ et les pavages (multiplicatifs) de nombres développés par Conway et Coxeter \cite{CC1,CC2} au cours des années 70. Bien que les pavages (multiplicatifs) de Conway et Coxeter soient bien connus, il existe certaines variantes de ceux-ci qui sont moins connues, par exemple les pavages additifs développés par G.C. Shephard \cite{S} ou les 2-frises étudiées par Morier-Genoud, Ovsienko et Tabachnikov \cite{MOT,SMG}. Ce qui suit est une introduction aux pavages additifs de Shephard comportant deux nouveaux théorèmes (2 et 5) et une conjecture.

\section{Pavages multiplicatifs}

Un \emph{pavage multiplicatif} est un assemblage de nombres arrang\'es en lignes d\'ecal\'ees comme les briques d'un mur de telle sorte que la ligne du haut et la ligne du bas sont compos\'ees uniquement de 1 et chaque losange
\begin{displaymath}
\begin{array}{ccc}
&a&\\
b&&c\\
&d&
\end{array}
\end{displaymath}
satisfait l'\'equation suivante:
\begin{equation}\label{uni}
\nonumber bc-ad=1.
\end{equation}

Cette équation est appel\'ee \emph{r\`egle unimodulaire}. Voici un exemple de pavage multiplicatif:
\begin{displaymath}
\begin{array}{ccccccccccccccccccc}
\ldots&1&&1&&1&&\textbf{1}&&1&&\textit{1}&&1&&1&&1&\ldots\\
&\ldots&1&&2&&2&&\textbf{2}&&2&&\textit{1}&&5&&1&\ldots&\\
\ldots&4&&1&&3&&3&&\textbf{3}&&1&&\textit{4}&&4&&1&\ldots\\
&\ldots&3&&1&&4&&\textbf{4}&&1&&3&&\textit{3}&&3&\ldots&\\
\ldots&2&&2&&1&&\textbf{5}&&1&&2&&2&&\textit{2}&&2&\ldots\\
&\ldots&1&&1&&1&&\textbf{1}&&1&&1&&1&&\textit{1}&\ldots
\end{array}
\end{displaymath}

\noindent\textbf{Définitions.}  Étant donné un pavage à $n$ lignes, une \emph{tranche} est un ensemble de $n$ éléments adjacents, où chaque élément provient d'une ligne différente (en gras dans l'exemple précédent). Un \emph{couple} est une paire d'éléments adjacents sur une tranche. Une \emph{diagonale} est une tranche à travers laquelle il est possible de tracer une ligne droite (en italique dans l'exemple précédent).

Dans ce document, nous utiliserons les notations suivantes afin de spécifier la position d'un nombre appartenant à un pavage:
\begin{displaymath}
\begin{array}{cccccccccccc}
\ldots&1&&1&&1&&1&&1&\ldots&$ligne 1$\\
&\ldots&m_{1,3}&&m_{2,4}&&m_{3,5}&&m_{4,6}&\ldots&&$ligne 2$\\
&&\ldots&m_{1,4}&&m_{2,5}&&m_{3,6}&\ldots&&&$ligne 3$\\
&&&\ldots&m_{1,5}&&m_{2,6}&\ldots&&&&$ligne 4$\\
&&&&\ldots&m_{1,6}&\ldots&&&&&$ligne 5$\\
&&&&&\ldots&&&&&&\vdots
\end{array}
\end{displaymath}
Note : Dans un pavage multiplicatif à $n$ lignes, la ligne 1 et la ligne $n$ sont compos\'ees uniquement de 1.

\subsection{Propri\'et\'es}

Les huit propriétés des pavages multiplicatifs qui suivent serviront de points de comparaisons avec les propriétés des pavages additifs. Lors de nos recherches nous avons remarqué que pour chaque propriété du cas multiplicatif il existe une propriété analogue pour le cas additif. Il est à noter que toutes les propriétés qui suivent proviennent des travaux de Conway et Coxeter \cite{CC1,CC2}, à l'exception de la propriété 3 qui a été tirée de l'ouvrage de Fraser Martineau et Lavertu \cite{FML}.

\begin{p1}
Si une diagonale du pavage est connue, il est possible de construire le reste du pavage et cette contruction est unique.
\end{p1}
\begin{p1} 
Tout pavage multiplicatif \`a $n$ lignes admet au moins $2$ isom\'etries: une translation de $n+1$ positions vers la droite (ou gauche) et une transvection compos\'ee d'une r\'eflexion par rapport à l'axe horizontal suivie d'une translation de $\frac{n+1}{2}$ positions vers la gauche (ou droite).
\end{p1}
\begin{p1}
Si un pavage multiplicatif contient une tranche composée uniquement de $1$, alors le pavage sera compos\'e uniquement d'entiers positifs.
\end{p1}
\begin{p1}
Les nombres sur la deuxi\`eme ligne d'un pavage multiplicatif \`a $n$ lignes correspondent aux nombres de triangles adjacents aux sommets d'un polygone triangul\'e (au sens de Conway et Coxeter \emph{\cite{CC1}}) \`a $n+1$ sommets.
\end{p1}
\begin{p1}
Le pavage sera compos\'e uniquement d'entiers positifs si et seulement s'il existe $k\in \mathbb{N}$ tel que pour tout $p\in \mathbb{N}$:
\begin{displaymath}
m_{k,p}~|~m_{k,p+1}+m_{k,p-1}.
\end{displaymath}
\end{p1}
\begin{p1}
\begin{displaymath}
\begin{array}{c|cccccccc|}
&m_{k,k+2}&1&0&\ldots&&&&\\
&1&m_{k+1,k+3}&1&0&\ldots&&&\\
&0&1&\ddots&&&&&\\
m_{k,p}=&\vdots&0&&&&&\vdots&\\
&&\vdots&&&&&0&\vdots\\
&&&&&&\ddots&1&0\\
&&&&\ldots&0&1&m_{p-3,p-1}&1\\
&&&&&\ldots&0&1&m_{p-2,p}
\end{array}
\end{displaymath}
\end{p1}
\begin{p1}
Pour tout $k\in \mathbb{N}$:
$$
m_{p,p+2}=\frac{m_{k,p}+m_{k,p+2}}{m_{k,p+1}}\cdot
$$
\end{p1}
\begin{p1}
Il est possible d'obtenir un pavage à $n+1$ lignes à partir d'un pavage à $n$ lignes en remplaçant 4 éléments consécutifs (\ldots,a,b,c,d,\ldots) de la deuxième ligne par (\ldots,a,b+1,1,c+1,d,\ldots).
\end{p1}

\section{Pavages additifs}

Un \emph{pavage additif} est un assemblage de nombres arrangés en lignes d\'ecal\'ees comme les briques d'un mur de telle sorte que la ligne du haut et la ligne du bas sont compos\'ees uniquement de 0 et chaque losange
\begin{displaymath}
\begin{array}{ccc}
&a&\\
b&&c\\
&d&
\end{array}
\end{displaymath}
satisfait l'\'equation suivante:
\begin{equation}\label{unimod}
(b+c)-(a+d)=1.
\end{equation}

\begin{exemple}
\begin{displaymath}
\begin{array}{ccccccccccccccc}
\ldots&0&&0&&0&&0&&0&&0&&0&\ldots\\
&\ldots&0&&1&&2&&3&&0&&1&\ldots&\\
\ldots&2&&0&&2&&4&&2&&0&&2&\ldots\\
&\ldots&1&&0&&3&&2&&1&&0&\ldots&\\
\ldots&0&&0&&0&&0&&0&&0&&0&\ldots
\end{array}
\end{displaymath}
\end{exemple}

\begin{exemple}
\begin{displaymath}
\begin{array}{ccccccccccccccc}
\ldots&0&&0&&0&&0&&0&&0&&0&\ldots\\
&\ldots&0&&3.5&&1.5&&1&&0&&3.5&\ldots&\\
\ldots&0&&2.5&&4&&1.5&&0&&2.5&&4&\ldots\\
&\ldots&1.5&&2&&3&&-0.5&&1.5&&2&\ldots&\\
\ldots&0&&0&&0&&0&&0&&0&&0&\ldots
\end{array}
\end{displaymath}
\end{exemple}
Pour la suite du document, nous utiliserons les notations suivantes pour situer les éléments dans un pavage additif:
\begin{displaymath}
\begin{array}{cccccccccccc}
\ldots&0&&0&&0&&0&&0&\ldots&$ligne 0$\\
&\ldots&a_{1}&&a_{2}&&a_{3}&&a_{4}&\ldots&&$ligne 1$\\
&&\ldots&a_{12}&&a_{23}&&a_{34}&\ldots&&&$ligne 2$\\
&&&\ldots&a_{123}&&a_{234}&\ldots&&&&$ligne 3$\\
&&&&\ldots&a_{1234}&\ldots&&&&&$ligne 4$\\
&&&&&\ldots&&&&&&\vdots
\end{array}
\end{displaymath}
Note : Dans un pavage additif à $n+1$ lignes, la ligne 0 et la ligne $n$ sont compos\'ees uniquement de 0. Dans le contexte des pavages additifs, la ligne $k$ sera aussi appelée la $k^e$ ligne. 

\subsection{Propri\'et\'es du pavage additif}
Cette section est l'essence du document, nous y présenterons les principaux théorèmes connus à propos des pavages additifs et nous fournirons des preuves à ceux-ci.

Le théorème suivant (tiré de \cite{S}) est l'analogue de la propriété 6 des pavages multiplicatifs.
\begin{theo}\label{somme}
Nous avons
\begin{displaymath}
a_{r\ldots s}=\left(\sum_{i=r}^{s}a_i\right) - t_{s-r},
\end{displaymath}
o\`u
\begin{displaymath}
t_{n}=\frac{n(n+1)}{2},
\end{displaymath}
c'est-à-dire que $t_{n}$ est le n$^{e}$ nombre triangulaire.
\end{theo}
\begin{proof}[Démonstration]

Nous procéderons par récurrence.

\noindent Pour $s-r=1$:
\begin{displaymath}
a_{rs}=a_{r}+a_{s}-1=\left(\displaystyle\sum_{i=r}^{s}a_i\right) - t_{1}.
\end{displaymath}

\noindent Pour $s-r=2$:
$$
\begin{array}{rl}
a_{rps}&=a_{rp}+a_{ps}-a_p-1\\[3mm]
&=(a_r+a_p-1)+(a_p+a_s-1)-a_p-1\\[3mm]
&=a_r+a_p+a_s-3\\[3mm]
&=\left(\displaystyle\sum_{i=r}^{s}a_i\right) - t_{2}.
\end{array}.
$$

\noindent Pour $s-r\geq3$: en appliquant l'équation (\ref{unimod}) sur le losange
$$
\begin{array}{ccc}
&a_{r+1\ldots s-1}&\\
a_{r\ldots s-1}&&a_{r+1...s}\\
&a_{r\ldots s}&
\end{array}
$$
nous obtenons:
$$
\begin{array}{rl}
a_{r\ldots s}&=a_{r\ldots s-1}+a_{r+1\ldots s}-a_{r+1\ldots s-1}-1\\ [3mm]
&=\left(\displaystyle\sum_{i=r}^{s-1}a_i+\displaystyle\sum_{i=r+1}^{s}a_i-\displaystyle\sum_{i=r+1}^{s-1}a_i\right)-2t_{s-r-1}+t_{s-r-2}-1\\ [5mm]
&=\left(\displaystyle\sum_{i=r}^{s}a_i\right)-\displaystyle\frac{2(s-r-1)(s-r)}{2}+\frac{(s-r-2)(s-r-1)}{2}-1\\ [5mm]
&=\left(\displaystyle\sum_{i=r}^{s}a_i\right)-\displaystyle\frac{(s^2-2sr+s+r^2-r)}{2}\\ [5mm]
&=\left(\displaystyle\sum_{i=r}^{s}a_i\right)-\displaystyle\frac{(s-r)(s-r+1)}{2}\\ [5mm]
&=\left(\displaystyle\sum_{i=r}^{s}a_i\right)-t_{s-r}.
\hfill\qedhere
\end{array}
$$
\end{proof}
\noindent Ce théorème apporte trois corollaires qui nous permettent de caractériser les pavages additifs et qui s'avéreront très importants pour démontrer certains théorèmes. Le premier corollaire est une condition suffisante pour que le pavage soit composé uniquement d'entiers.
\begin{coro}\label{coro1}
Le pavage est uniquement composé d'entiers si et seulement si les éléments de la première ligne sont tous des entiers.
\end{coro}
Le second corollaire est une condition suffisante pour qu'un pavage additif en soit un, c'est-à-dire que la $n^e$ ligne soit composée uniquement de 0.
\begin{coro}\label{coro2}
La $n^e$ ligne est une ligne de 0 si et seulement si la somme de n'importe quelle suite de $n$ nombres consécutifs sur la première ligne est égale à $t_{n-1}$.
\end{coro}
Le dernier corollaire est une condition suffisante pour que le pavage additif soit composé uniquement de nombres positifs.
\begin{coro}\label{coro3}
Si pour tous $r,s \in \mathbb{Z}$ avec $s \leq n$, nous avons
$$\sum_{i=r}^{r+s}a_i \geq t_s$$
alors le pavage additif est composé uniquement de nombres positifs.
\end{coro}
Le théorème précédent nous a montré une relation permettant de calculer un élément quelconque du pavage à partir des éléments sur la première ligne. Maintenant, il serait intéressant de voir s'il existe une relation permettant de calculer les éléments de la première ligne à partir d'éléments quelconques. En fait, si un couple d'éléments est connu, il existe une telle relation. Celle-ci est présentée dans le théorème suivant qui se veut l'analogue de la propriété 7 des pavages multiplicatifs.
\begin{theo}\label{mien}
Pour tous $r,s$, avec $s>r$, nous avons:
\begin{enumerate}
\item[\emph{(a)}] $a_s =  a_{r\ldots s}-a_{r\ldots s-1}+(s-r)$, 
\item[\emph{(b)}] $a_r =  a_{r\ldots s}-a_{r+1\ldots s}+(s-r)$.
\end{enumerate}
\end{theo}

\begin{proof}[Démonstration]
(a) En vertu du théorème \ref{somme}, nous avons:
$$
a_{r\ldots s}-a_{r\ldots s-1}+(s-r)=\left(\sum_{i=r}^{s}a_i\right)-\left(\sum_{i=r}^{s-1}a_i\right)-t_{s-r}+t_{s-r-1}+(s-r)=a_s
$$

(b) La preuve est similaire et laissée au lecteur.
\end{proof}

Le théorème suivant, tiré de \cite{S},  est l'analogue de la propriété 2 des pavages multiplicatifs.
\begin{theo}
Tout pavage additif à $n+1$ lignes admet comme isométrie une translation de $n$ positions vers la droite (ou la gauche).
\end{theo}
\begin{proof}[Démonstration]
En vertu du théorème \ref{somme}, nous savons que:
$$
\left(\sum_{i=k}^{k+n-1}a_i\right)-t_{n-1}=0=\left(\sum_{i=k+1}^{k+n}a_i\right)-t_{n-1}.
$$
En simplifiant, on obtient
$$
a_k=a_{k+n}.
$$
Ensuite,
$$
\begin{array}{rl}
a_{r\ldots s}&=\left(\displaystyle\sum_{i=r}^{s}a_i\right)-t_{s-r}\\[3mm]
&=\left(\displaystyle\sum_{i=r+n}^{s+n}a_i\right)-t_{s-r}\\[3mm]
&=a_{r+n\ldots s+n}.
\end{array}
$$
\end{proof}
Le théorème suivant, également tiré de \cite{S}, est l'analogue de la propriété 3 des pavages multiplicatifs.
\begin{theo}
Si un pavage additif à $n+1$ lignes admet une tranche de 0, alors le pavage est composé uniquement d'entiers positifs.
\end{theo}
\begin{proof}[Démonstration]
En appliquant le théorème \ref{mien} sur tous les couples de zéro qui composent la tranche, nous pouvons calculer $n$ nombres consécutifs sur la première ligne. De cette manière, nous découvrons que la première ligne est en fait  une permutation des nombres de 1 à $n$. Comme toutes les permutations de 1 à $n$ répondent aux trois corollaires du théorème \ref{somme}, le pavage est composé uniquement d'entiers positifs.
\end{proof}
Le théorème suivant est l'analogue de la propriété 8 des pavages multiplica\-tifs.
\begin{theo} \label{mien 2}
Il est possible d'obtenir un pavage à $n+1$ lignes à partir d'un pavage à $n$ lignes en insérant le nombre $n-1$ n'importe où dans la première ligne du pavage. De plus, si le pavage à $n$ lignes est formé d'entiers positifs, alors il en est de même pour le pavage à $n+1$ lignes.
\end{theo}
\begin{proof}[Démonstration]
Pour démontrer la première partie de l'énoncé, il suffit de vérifier que les hypothèses du corollaire 2 sont satisfaites pour le nouvel arrangement à $n+1$ lignes. Supposons un pavage additif à $n$ lignes. Nous savons que:
$$
\begin{array}{rrl}
&\displaystyle\left(\sum_{i=k}^{k+n-2}a_i\right)&=t_{n-2}
\end{array}
$$
pour tout $k$. Alors:
$$
\displaystyle\left(\sum_{i=k}^{k+n-2}a_i\right)+n-1 = t_{n-2} + n-1 = t_{n-1}.
$$
de sorte que les hypothèses du corollaire 2 sont vérifiées pour le nouvel arrangement. On obtient donc un pavage à $n+1$ lignes.

Pour ce qui est de la seconde partie, il suffit de montrer que les hypothèses des corollaires 1 et 3 sont satisfaites pour le nouveau pavage à $n+1$ lignes. Pour montrer que les hypothèses du corollaire 3 sont satisfaites, il suffit de le montrer pour les sous-sommes qui incluent le nombre $n-1$ nouvellement ajouté à la première ligne. En vertu du théorème \ref{somme}, nous savons que:
$$
\left(\sum_{i=r}^{r+s}a_i\right)\geq t_s,
$$
d'où
$$
\left(\displaystyle\sum_{i=r}^{r+s}a_i\right)+(n-1)\geq t_s + (n-1)\geq t_s + (s+1)=t_{s+1}
$$
car $s$ est au maximum égal à $n-2$. Les hypothèses du corollaire 3 sont donc vérifiées pour le pavage à $n+1$ lignes. Finalement, puisque $n-1$ est un entier, l'hypothèse du corollaire 1 est vérifiée.
\end{proof}
\noindent Note: Cette méthode fonctionne sur tous les pavages mais il ne permet pas de tous les obtenir.

\begin{exemple}
Considérons le pavage original suivant:
$$
\begin{array}{ccccccccccccccc}
0&&0&&0&&0&&0&&0&&0&&0\\
&2&&2&&2&&2&&2&&2&&2&\\
3&&3&&3&&3&&3&&3&&3&&3\\
&3&&3&&3&&3&&3&&3&&3&\\
2&&2&&2&&2&&2&&2&&2&&2\\
&0&&0&&0&&0&&0&&0&&0&
\end{array}
$$
En insérant un 5 dans la première rangée nous obtenons:
$$
\begin{array}{ccccccccccccccc}
0&&0&&0&&0&&0&&0&&0&&0\\
&2&&2&&5&&2&&2&&2&&2&\\
3&&3&&6&&6&&3&&3&&3&&3\\
&3&&6&&6&&6&&3&&3&&3&\\
2&&5&&5&&5&&5&&2&&2&&5\\
&3&&3&&3&&3&&3&&0&&3&\\
0&&0&&0&&0&&0&&0&&0&&0
\end{array}
$$
Toutefois, il est évident qu'il est impossible d'obtenir le pavage original de cette façon, puisqu'il contient 5 lignes et qu'il n'y a aucun 4 sur la première ligne.
\end{exemple}

Le théorème qui suit provient de \cite{S}.
\begin{theo}
Supposons un parallélogramme ayant pour sommets: $$(a_{j\ldots k}, a_{j\ldots s}, a_{r\ldots s}, a_{r\ldots k})$$ 
où $r < j \leq k < s$. Alors:
$$
a_{j\ldots k}+a_{r\ldots s}=a_{j\ldots s}+a_{r\ldots k}-(j-r)(s-k),
$$
où $(j-r)(s-k)$ représente l'aire du parallélogramme si nous définissons l'aire des losanges comme étant égale à 1.
\end{theo}
\begin{proof}[Démonstration]
$$
\begin{array}{lll}
a_{r\ldots k}& + & a_{j\ldots s}-(j-r)(s-k)\\[5mm]
& = & \left(\left(\sum_{i=r}^{k}a_i\right)-t_{k-r}\right)+\left(\left(\sum_{i=j}^{s}a_i\right)-t_{s-j}\right)-(j-r)(s-k)\\  [5mm]
&= & \left(\sum_{i=r}^{k}a_i\right)+\left(\sum_{i=j}^{s}a_i\right)-\frac{(k-r)(k-r+1)}{2}-\frac{(s-j)(s-j+1)}{2}\\ [5mm]
& & -(j-r)(s-k)\\ [5mm]
&= & \Big(\sum_{i=r}^{s}a_i\Big)+\left(\sum_{i=j}^{k}a_i\right)-\frac{(k^2-2kr+k-r+r^2)}{2} -\frac{(s^2-2sj+s-j+j^2)}{2}\\ [5mm]
& &-\frac{2(js-jk-rs+rk)}{2}\\ [5mm]
&=& \Big(\sum_{i=r}^{s}a_i\Big)+\left(\sum_{i=j}^{k}a_i\right)+\frac{(-k^2-k+r-r^2-s^2-s+j-j^2+2jk+2rs)}{2}\\ [5mm]
&=& \Big(\sum_{i=r}^{s}a_i\Big)+\left(\sum_{i=j}^{k}a_i\right)-\frac{(k^2-2jk+k-j+j^2)}{2}-\frac{(s^2-2sr+s-r+r^2)}{2} \\[5mm]
&=& \Big(\sum_{i=r}^{s}a_i\Big)+\left(\sum_{i=j}^{k}a_i\right)-\frac{(k-j)(k-j+1)}{2}-\frac{(s-r)(s-r+1)}{2}\\ [5mm]
&=& \Big(\sum_{i=r}^{s}a_i\Big)+\left(\sum_{i=j}^{k}a_i\right)-t_{k-j}-t_{s-r}\\ [5mm]
&=& a_{j\ldots k}+a_{r\ldots s}
\end{array}
$$
\end{proof}
\section{Pavages constants}
Un \emph{pavage constant} est un pavage (additif ou multiplicatif) pour lequel sur chaque ligne tous les éléments sont égaux.\\

\noindent \textbf{Exemple additif:}
$$
\begin{array}{ccccccccccc}
&0&&0&&0&&0&&0&\\
&&2&&2&&2&&2&&\\
\ldots&3&&3&&3&&3&&3&\ldots\\
&&3&&3&&3&&3&&\\
&2&&2&&2&&2&&2&\\
&&0&&0&&0&&0&&
\end{array}
$$

\noindent \textbf{Exemple multiplicatif:}
$$
\begin{array}{ccccccccccc}
&1&&1&&1&&1&&1&\\
&&x&&x&&x&&x&&\\
\ldots&y&&y&&y&&y&&y&\ldots\\
&&x&&x&&x&&x&&\\
&1&&1&&1&&1&&1&
\end{array}
$$
où $x=\sqrt{3}$ et $y=2$.\\
Une interprétation géométrique peut être associée à chacun des deux types de pavages.

\begin{figure}[h]
\centering
\includegraphics[width=4cm]{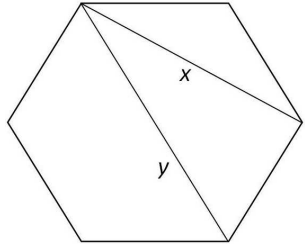}
\caption{Hexagone régulier avec des côtés de longueur 1.}
\label{fig:hexagone}
\end{figure}

\medskip
\noindent \textbf {Le cas multiplicatif:}\\
Pour le cas multiplicatif, il a été démontré dans \cite{CC1,CC2} que la valeur des éléments sur une ligne d'un pavage multiplicatif constant à $n$ lignes représente la longueur des diagonales d'un polygone régulier à $n+1$ côtés de longueur 1. L'exemple précédent illustre les deux longueurs possible pour les diagonales d'un hexagone régulier avec des côtés de longueur 1 (voir Figure~\ref{fig:hexagone}).

\begin{figure}[h]
\centering
		\includegraphics[width=5cm]{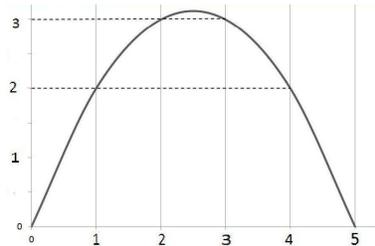}
		\caption{$f(x)=\frac{5x-x^2}{2}$}
	\label{fig:parabole}
\end{figure}

\noindent \textbf {Le cas additif:}\\
 Il est spécifié dans \cite{S} que tout pavage additif constant est directement relié à une parabole.  Par exemple, le pavage 
$$
\begin{array}{ccccccccccc}
&0&&0&&0&&0&&0&\\
&&2&&2&&2&&2&&\\
\ldots&3&&3&&3&&3&&3&\ldots\\
&&3&&3&&3&&3&&\\
&2&&2&&2&&2&&2&\\
&&0&&0&&0&&0&&
\end{array}
$$
\noindent peut être relié à la parabole présentée à la Figure~\ref{fig:parabole}.

Cet exemple met en évidence que la valeur sur la $n$-ième ligne du pavage est égale à la valeur en $y$ de la parabole pour $x=n$. Cette relation est applicable pour tous les pavages additifs constants. L'équation de la parabole correspondant à un pavage additif constant à $n$ lignes est donnée par 
$$
f(x)=\frac{nx-x^2}{2}\cdot
$$
\section{Conclusion et conjecture}
Pour le cas des pavages additifs constants, la relation avec les paraboles est intéressante car elle est le seul lien connu avec un autre sujet mathématique et peut être une piste pour découvrir une interprétation géométrique aux pavages additifs en général. Lors de nos recherches à propos des pavages, nous nous sommes demandé s'il était possible de trouver une nouvelle interprétation géométrique aux lignes centrales des pavages multiplicatifs à plus de quatre lignes. Nous n'avons pu en trouver une, mais nous en sommes arrivés à une conjecture que nous n'avons pas pu démontrer. Nous avons découvert que si nous prenons une ligne autre que la première, deuxième, la dernière ou l'avant-dernière ligne d'un pavage à cinq lignes ou plus et que nous l'associons à une ligne de 1, nous obtenons un pavage infini (n'ayant pas de ligne de 1 à la fin) composé uniquement d'entiers positifs. Par exemple, si nous prenons la troisième ligne du pavage suivant:
\begin{displaymath}
\begin{array}{ccccccccccccccccccc}
\ldots&1&&1&&1&&1&&1&&1&&1&&1&&1&\ldots\\
&\ldots&1&&2&&2&&2&&2&&1&&5&&1&\ldots&\\
\ldots&4&&1&&3&&3&&3&&1&&4&&4&&1&\ldots\\
&\ldots&3&&1&&4&&4&&1&&3&&3&&3&\ldots&\\
\ldots&2&&2&&1&&5&&1&&2&&2&&2&&2&\ldots\\
&\ldots&1&&1&&1&&1&&1&&1&&1&&1&\ldots
\end{array}
\end{displaymath}
et que nous adjoignons cette ligne à une ligne de 1 nous obtenons:
\begin{displaymath}
\begin{array}{ccccccccccccccccc}
\ldots&1&&1&&1&&1&&1&&1&&1&&1&\ldots\\
&\ldots&4&&1&&3&&3&&3&&1&&4&\ldots&\\
\ldots&15&&3&&2&&8&&8&&2&&3&&15&\ldots\\
&\ldots&11&&5&&5&&21&&5&&5&&11&\ldots&\\
\ldots&8&&18&&12&&13&&13&&12&&18&&8&\ldots\\
&&13&&43&&31&&8&&31&&43&&13&&\\
&\vdots&&\vdots&&\vdots&&\vdots&&\vdots&&\vdots&&\vdots&&\vdots&\\
\end{array}
\end{displaymath}
Ce genre de pavages pourrait faire l'objet de recherches plus approfondies pour de futurs projets. 
\section*{Remerciements}
Cette recherche a été rendue possible grâce à une bourse de recherche de premier cycle du CRSNG et grâce à David Smith de l'Université Bishop's qui m'a guidé tout au long de mes recherches.

\renewcommand\refname{Références}
\end{document}